\documentclass[A4paper]{article}
\usepackage{amsmath}
\usepackage{amssymb}
\usepackage{amsthm}
\usepackage{graphicx,psfrag,epsf}
\usepackage{enumerate}
\usepackage{natbib}
\usepackage{url} % not crucial - just used below for the URL 
\usepackage{fullpage}

\newcommand{\blind}{1}

\renewcommand\footnotemark{}

\def\T{{ \mathrm{\scriptscriptstyle T} }}

\newcommand{\bt}{\mathbf{t}}
\newcommand{\bd}{\mathbf{d}}
\newcommand{\bD}{\mathbf{D}}
\newcommand{\by}{\mathbf{y}}
\newcommand{\bv}{\mathbf{v}}

\newcommand{\bmu}{\boldsymbol{\mu}}

\newcommand{\btheta}{\boldsymbol{\theta}}

\newcommand{\bone}{\mathbf{1}}

\newtheorem{theorem}{Theorem}[section]

\newtheorem{lemma}[theorem]{Lemma}
\newtheorem{proposition}[theorem]{Proposition}
\begin{document}

\def\spacingset#1{\renewcommand{\baselinestretch}%
{#1}\small\normalsize} \spacingset{1}

%%%%%%%%%%%%%%%%%%%%%%%%%%%%%%%%%%%%%%%%%%%%%%%%%%%%%%%%%%%%%%%%%%%%%%%%%%%%%%

\if1\blind
{
  \title{\bf \Large Properties of Fisher information gain for Bayesian design of experiments
  \vspace{0cm}
  }
  \author{\large Antony M. Overstall\hspace{.2cm}\\
    \large Southampton Statistical Sciences Research Institute,\\ \large University of Southampton,\\ \large Southampton SO17 1BJ UK \\ \large (A.M.Overstall@southampton.ac.uk)}
    \date{\vspace{-1.2cm}}
  \maketitle
} \fi

\if0\blind
{
  \bigskip
  \bigskip
  \bigskip
  \begin{center}
    {\LARGE\bf Title}
\end{center}
  \medskip
} \fi

\bigskip
\begin{abstract}
The Bayesian decision-theoretic approach to design of experiments involves specifying a design (values of all controllable variables) to maximise the expected utility function (expectation with respect to the distribution of responses and parameters).  For most common utility functions, the expected utility is rarely available in closed form and requires a computationally expensive approximation which then needs to be maximised over the space of all possible designs. This hinders practical use of the Bayesian approach to find experimental designs. However, recently, a new utility called Fisher information gain has been proposed. The resulting expected Fisher information gain reduces to the prior expectation of the trace of the Fisher information matrix. Since the Fisher information is often available in closed form, this significantly simplifies approximation and subsequent identification of optimal designs. In this paper, it is shown that for exponential family models, maximising the expected Fisher information gain is equivalent to maximising an alternative objective function over a reduced-dimension space, simplifying even further the identification of optimal designs. However, if this function does not have enough global maxima, then designs that maximise the expected Fisher information gain lead to non-identifiablility.
\end{abstract}

\noindent%

\section{Introduction}
\label{sec:INTRO}

Suppose an experiment is to be conducted to learn the relationship between a series of controllable variables and a measurable response for some physical system. The experiment consists of a number of runs where each run involves specification of a series of controllable variables and subsequent observation of a response. Typically, on completion of the experiment, it is assumed that the observed responses are realisations of a random variable with a fully specified probability distribution apart from a vector of unknown parameters. That is, a statistical model is assumed for the responses. The observed responses are used to estimate the parameters thus allowing estimation of the relationship between controllable variables and response.

A question of interest is how should the controllable variables for all runs (i.e. the design) be specified to facilitate as precise as possible estimation of the unknown parameters. The Bayesian decision-theoretic approach \citep[see, for example,][]{CV1995} is to select the design to maximise the expectation of a utility function. Expectation is taken with respect to the joint distribution of all unknown quantities, i.e. the responses and the parameters, with this joint distribution following from the specification of the statistical model and a prior distribution for the parameters. The advantages of this approach are as follows. (1) The choice of utility allows bespoke experimental aims to be incorporated (for example, experimental aims of prediction and/or model selection can be considered instead of estimation although these are not the focus of this paper). (2) By taking expectation with respect to unknown quantities, all known sources of uncertainty are incorporated. (3) The framework fits into the iterative nature in which knowledge is accumulated in science, i.e. the distribution of unknown quantities for the current iteration is based on the results of the previous iteration. With regards to point (3) above, the approach can even be extended to account for the fact that the expected utility will be maximised in subsequent iterations \citep[see, for example,][]{huanmarzouk_2016}, i.e. so-called Bayesian sequential or adaptive design. However this extension is not considered in this paper. 

%First, note that the use of Bayesian methods to determine a design does not restrict the post-experiment analysis to be Bayesian \citep[see, for example,][Chapter 18]{ADT2007}. Second, note that a design that maximises an expected utility is not necessarily unique.

The choice of utility allows different ways of specifying estimation precision. Ideally, the choice of utility should be uniquely tailored to the experiment. However, a more pragmatic choice is often made from certain default utility functions. Two common default utilities for estimation precision are Shannon information gain \citep{lindley1956} and negative squared error loss \citep[see, for example,][pages 77-79]{robert2007}. In either case, typically, the utility function depends on the responses through the posterior distribution of the parameters. Even for relatively simple statistical models, the posterior distribution is not of a known form meaning that the utility and expected utility are not available in closed form and require approximation. 

Recently, \cite{walker2016} instead proposed the Fisher information gain utility. This is the difference between the squared Euclidean length of the gradient of the log posterior and log prior densities of the parameters. The appealing feature of Fisher information gain is that the resulting expected utility is the prior expectation of the trace of the Fisher information matrix. For many models, the Fisher information matrix is available in closed form, significantly simplifying approximation of the expected utility and subsequent identification of an optimal design.

However, in this paper it is shown that, under a broad class of models, a design that maximises the expected Fisher information gain can lead to a non-identifiable model. Under the classical approach to statistical inference, this means the parameters are not uniquely estimable \citep[see, for example, ][]{catchpolemorgan1997}. Under the Bayesian approach, it means that the posterior and prior distributions for a subset of parameters (conditional on the complement of the subset) are identical \citep[see, for example, ][]{gelfandsahu_1999}. These are undesirable consequences for a default utility function.

The paper is organised as follows. Section~\ref{sec:BACK} provides a background to Bayesian decision-theoretic design of experiments, introduces the Fisher information gain utility, describes the exponential family of models and discusses the notion of non-identifiability. Section~\ref{sec:METH} explores the consequences of designing experiments under Fisher information gain with examples provided in Section~\ref{sec:EXAMPLES}.

\section{Background} \label{sec:BACK}

\subsection{Bayesian decision-theoretic design of experiments} \label{sec:BDOE}

Suppose there are $k$ controllable variables denoted $\bd = \left(d_1, \dots, d_k \right)^\T \in \mathcal{D}$, where $\mathcal{D}$ denotes a set of possible values for the controllable variables, and let $y$ denote the response.  Suppose the experiment consists of a fixed number $n$ of runs. For the $i$th run, the controllable variables are denoted $\bd_i = \left(d_{i1}, \dots, d_{ik} \right)^{\T}$ and the response by $y_i$, for $i=1,\dots,n$. Let $\by = \left(y_1,\dots,y_n\right)^T$ denote the $n \times 1$ vector of responses and $\bD = \left(\bd_1^\T,\dots, \bd_n^\T \right)^\T \in \mathcal{D}^n$ denote the $nk \times 1$ vector of controllable variables for all runs (termed the design). The vectors in the set $\left\{ \bd_1, \dots, \bd_n \right\}$ are referred to as designs points and the distinct vectors as support points.

Throughout it is assumed that the elements of $\by$ are independent and $y_i \sim \mathrm{F}(\btheta, \bd_i)$ where $\mathrm{F}(\btheta, \bd)$ denotes a probability distribution depending on a $p \times 1$ vector of unknown parameters $\btheta = \left(\theta_1, \dots, \theta_p \right)^T \in \Theta$ and controllable variables $\bd$, with $\Theta$ denoting a $p$-dimensional parameter space. It is assumed that $p \le n$.

Bayesian decision-theoretic design of experiments \citep[see, for example, ][]{CV1995} proceeds by specifying a utility function denoted $u(\btheta, \by, \bD)$. The utility gives the precision in estimating $\btheta$ given responses $\by$ obtained via design $\bD$, where different choices of utility allow different ways of measuring precision. Typically the utility depends on the responses through the posterior distribution of $\btheta$ given by Bayes' theorem $\pi(\btheta \vert \by, \bD) \propto \pi(\by \vert \btheta, \bD) \pi(\btheta)$, where $\pi(\btheta)$ is the probability density function (pdf) of the prior distribution of $\btheta$ (assumed independent of the design $\bD$) and $\pi(\by \vert \btheta, \bD) = \prod_{i=1}^n \pi(y_i \vert \btheta, \bd_i)$ is the likelihood function with $\pi(y_i \vert \btheta, \bd_i)$ the pdf or mass function of $F(\btheta, \bd_i)$.

Let
$$U(\bD) = \mathrm{E}_{\by, \btheta \vert \bD} \left[ u(\btheta, \by, \bD) \right],$$
be the expected utility where expectation is with respect to the joint distribution of responses $\by$ and parameters $\btheta$. Assuming that $U(\bD)$ exists, an optimal Bayesian design maximises $U(\bD)$ over the space of all designs $\mathcal{D}^n$, where this design is not necessarily unique. 

It should be noted that the use of Bayesian methods to determine a design does not restrict the post-experiment analysis to be Bayesian \citep[see, for example,][Chapter 18]{ADT2007}.

Ideally, the utility should be chosen to represent the experimental aim. Instead, often a more pragmatic approach is taken where the utility is chosen from certain default utility functions. These utilities aim to represent standard aims such as precision of point estimation or gain in information from prior to posterior. Common default choices for the utility are Shannon information gain (SIG) $u_{SIG}(\btheta, \by, \bD) = \log \pi(\btheta \vert \by, \bD) - \log \pi(\btheta)$ and negative squared error loss (NSEL) $u_{NSEL}(\btheta, \by, \bD) = - \Vert \btheta - \mathrm{E}_{\btheta \vert \by, \bD} \left( \btheta \right) \Vert$, where $\Vert \bv \Vert = \bv^\T \bv$ is the squared Euclidean length of vector $\bv$. A SIG-optimal design that maximises the expected SIG utility equivalently maximises the expected Kullback-Liebler divergence between the posterior and prior distributions of $\btheta$, where expectation is with respect to the marginal distribution of $\by$. A NSEL-optimal design equivalently minimises the trace of the expected posterior variance matrix of $\btheta$.

The dependence of the utility function on the posterior distribution of $\btheta$, and the fact that this distribution is typically not of known form for non-trivial statistical models, is a significant hurdle to the practical implementation of Bayesian design of experiments. For example, under SIG and NSEL, the utility function is not available in closed form, not to mention the expected SIG and NSEL utilities. Even with the ongoing development of new methodology \citep[see, for example,][]{ryanetal2016, woodsetal2017} to approximate and maximise the expected utility, finding designs in practice is a computationally expensive task.

\subsection{Fisher information gain} \label{sec:FIG}

\cite{walker2016} recently proposed the Fisher information gain (FIG) utility
$$u_{FIG}(\btheta, \by, \bD) = \left\Vert \frac{\partial \log \pi(\btheta \vert \by, D)}{\partial \btheta}\right\Vert - \left\Vert \frac{\partial \log \pi(\btheta)}{\partial \btheta}\right\Vert,$$
the difference between the squared Euclidean length of the $p \times 1$ gradient vectors of the posterior and prior log pdfs (assuming that these gradients exist). The length of the gradient vector of the log pdf of a continuous distribution provides a measure of information provided by the distribution. The FIG utility gives the gain in this information provided by the posterior compared to the prior. A design that maximises the expected Fisher information gain $U_{FIG}(\bD) = \mathrm{E}_{\by, \btheta \vert \bD} \left[ u_{FIG}(\btheta, \by, \bD) \right]$ is termed FIG-optimal and denoted $\tilde{\bD}$. 

Under standard regularity conditions \citep[see, for example,][page 111]{shervish1995}, it can be shown \citep{walker2016} that the expected FIG utility can be written
$$
U_{FIG}(\bD) = \mathrm{E}_{\btheta} \left[ \mathrm{tr} \left\{ \mathcal{I}(\btheta , \bD) \right\} \right],
$$
where
$$\mathcal{I}(\btheta , \bD) = \mathrm{E}_{\by \vert \btheta, \bD} \left[ \frac{\partial \log \pi(\by \vert \btheta, D)}{\partial \btheta} \frac{\partial \log \pi(\by \vert \btheta, D)}{\partial \btheta^\T} \right]$$
is the $p \times p$ Fisher information matrix. That is, the expected FIG utility is the prior expectation of the trace of the Fisher information matrix where it is assumed that both the Fisher information, and the prior expectation of its trace, exists. 

The Fisher information matrix is available in closed form for many classes of model (for example, generalised linear models) for which the corresponding posterior distribution (under any prior distribution) is not of known form. Therefore, approximating the expected FIG utility is a significantly simpler task than, for example, approximating the expected SIG or NSEL utilities. This represents a major advantage in the use of FIG over other utility functions.

Furthermore, the expected FIG appears to parallel common objective functions for pseudo-Bayesian classical design criteria, e.g. D- and A-optimality \citep[see, for example,][Chapter 10]{ADT2007}, where in each case, the objective function is given as the prior expectation of a scalar function of the Fisher information matrix. For D-optimality the scalar function is the log-determinant and, for A-optimality, it is the negative trace of the inverse. Under these classical criteria, the Fisher information appears since its inverse is an asymptotic approximation to the variance of the maximum likelihood estimator of $\btheta$. To overcome dependence on $\btheta$, the scalar function of $\mathcal{I}(\btheta, D)$ is averaged with respect to a prior distribution on $\btheta$. The term ``pseudo-Bayesian" refers to the use of Bayesian machinery for a classical procedure. However, the objective functions for D- and A-optimality can also be derived as asymptotic approximations to the expected utility under the SIG and NSEL utilities, respectively \citep[see, for example,][]{CV1995}. As pointed out by \cite{walker2016}, the FIG utility is fully Bayesian and does not result from an asymptotic approximation. For the objective function for D- and A-optimality to be greater than $- \infty$, the inverse of the Fisher information must exist. There is no such requirement for the expected FIG.

\subsection{Exponential family models} \label{sec:MODELS}

We now present a broad class of models for which the results in this paper apply. It is assumed that the elements of $\by$ are distributed according to a member of the exponential family of distributions \citep[see, for example][]{kosmidis_etal_2020}. This family of distributions include the normal, Poisson, binomial and gamma. The pdf or mass function of $F(\btheta, \bd_i)$ is 
%$$
%\pi(y_i \vert \zeta_i, \gamma) = \exp \left( \frac{y_i \zeta_i - b(\zeta_i)}{\gamma} + c_2(\gamma, y_i) \right),
%$$
$$
\pi(y_i \vert \zeta_i, \gamma) = \exp \left( \frac{y_i \zeta_i - b(\zeta_i) - c_1(y_i)}{\gamma} - \frac{1}{2} a\left(-\frac{1}{\gamma}\right) + c_2(y_i) \right),
$$
for $i=1,\dots,n$, where $a(\cdot)$, $b(\cdot)$, $c_1(\cdot)$ and $c_2(\cdot)$ are sufficiently smooth specified functions, $\gamma>0$ is a common scale parameter and $\zeta_1,\dots,\zeta_n$ are natural parameters, depending on $\btheta$ and $\bd_1, \dots,\bd_n$, respectively. It follows that $\mathrm{E}_{y_i \vert \zeta_i, \gamma}(y_i) = \mathrm{E}_{y_i \vert \btheta, \bd_i, \gamma}(y_i) = b'(\zeta_i)$ \citep[see, for example][]{kosmidis_etal_2020}, where it is assumed that $\mathrm{E}_{y_i \vert \zeta_i, \gamma}(y_i) = \mathrm{E}_{y_i \vert \btheta, \bd_i, \gamma}(y_i) = \mu(\btheta, \bd_i)$: a differentiable mean function of parameters $\btheta$ and controllable variables $\bd_i$. Additionally, $\mathrm{var}_{y_i \vert \zeta_i, \gamma}(y_i) = \mathrm{var}_{y_i \vert \btheta, \bd_i, \gamma}(y_i) = \gamma b''(\zeta_i)$.

The distribution of $y_i$ is uniquely determined by $\mu(\btheta, \bd_i)$ through the natural parameter $\zeta_i$. This means that the difference in distribution of elements $y_i$ and $y_j$ is entirely determined by $\bd_i$ and $\bd_j$, for $i \ne j$. 

The decision to focus on exponential family models is due to 1) a large number of statistical models used to analyse experimental responses belong to this family; 2) the Fisher information matrix exists and is available in closed form for models in this family making the use of Fisher information gain particularly attractive; and 3) the existence of results to determine whether or not an exponential family model is identifiable (see Section~\ref{sec:IDENT}).

If the scale parameter $\gamma$ is known, the Fisher information is the $p \times p$ matrix
$$
\mathcal{I}(\btheta, \bD) = M(\btheta, \bD)^\T W(\btheta, \bD) M(\btheta, \bD) 
$$
\citep[see, for example][]{kosmidis_etal_2020} where $M(\btheta, \bD)$ is an $n \times p$ matrix with $ij$th element given by 
\begin{equation}
M_{ij}(\btheta, \bD) = \frac{\partial \mu(\btheta, \bd_i)}{\partial \theta_j},
\label{eqn:Mmat}
\end{equation}
and $W(\btheta, \bD)$ is an $n \times n$ diagonal matrix with $i$th diagonal element
$$W_{ii}(\btheta, \bD) = \frac{1}{\gamma b''(\zeta_i)} = \frac{1}{\mathrm{var}_{y_i \vert \zeta_i, \gamma} (y_i)} = \frac{1}{\mathrm{var}_{y_i \vert \btheta, \bd_i, \gamma}(y_i)},$$
for $i=1,\dots,n$. 

If $\gamma$ is unknown and of interest, an extra row and column can be included in the Fisher information matrix. However, it is arguably more realistic to allow $\gamma$ to be unknown but not of direct interest, i.e. it is a nuisance parameter. The correct approach here is to obtain the marginal likelihood by marginalising over $\gamma$ with respect to its prior distribution, i.e.
\begin{equation}
\pi(\by \vert \btheta, \bD) = \int_0^\infty \pi(\by \vert \btheta, \gamma, \bD) \pi(\gamma) \mathrm{d} \gamma.
\label{eqn:marglik}
\end{equation}
where $\pi(\by \vert \btheta, \gamma, \bD) = \prod_{i=1}^n \pi(y_i \vert \zeta_i, \gamma)$ and $\pi(\gamma)$ is the pdf of the prior distribution of $\gamma$ (assumed independent of $\btheta$). The Fisher information matrix can then be obtained from the marginal likelihood (\ref{eqn:marglik}). However, the exact form of the marginal likelihood depends on the form of $a(\cdot)$, $b(\cdot)$, $c_1(\cdot)$ and $c_2(\cdot)$. In Section~\ref{sec:METH}, we consider the special case of a normal distribution.

\subsection{Non-identifiability and parameter redundancy} \label{sec:IDENT}

\cite{catchpolemorgan1997} define identifiability to mean that there do not exist two parameter values (separated by at least $\delta$ under some distance measure) such that the probability distribution for the responses is identical under both parameter values. An obvious cause of non-identifiability is parameter redundancy, i.e. the statistical model can be equivalently re-expressed using a subset of the parameters. For the exponential family models of Section~\ref{sec:MODELS}, \cite{catchpolemorgan1997} show that the model is parameter redundant if and only if the $n \times p$ matrix $M(\btheta,\bD)$ is rank deficient, i.e. the number of linearly independent rows (or columns) of $M(\btheta,\bD)$ is less than $p$.

For classical inference, a consequence of parameter redundancy is that the parameters are not uniquely estimable since the likelihood can exhibit a ridge, i.e. a curve of parameter values that all maximise the likelihood function.

For Bayesian inference, the issue is more nuanced. Under a proper prior distribution, the posterior distribution for $\btheta$ will exist for a parameter redundant model. However, parameter redundancy does result in Bayesian non-identifiability \citep[see, for example, ][]{gelfandsahu_1999}. This is where $\btheta$ can be decomposed as $\btheta = \left(\btheta_1^\T, \btheta_2^\T \right)^\T$ and the posterior distribution of $\btheta_1$, conditional on $\btheta_2$, is equal to the prior distribution of $\btheta_1$, conditional on $\btheta_2$. This does not mean that no learning has taken place about $\btheta_1$, rather that learning comes indirectly through $\btheta_2$.

\section{Maximising expected Fisher information gain} \label{sec:METH}

In this section we show that the design points of a FIG-optimal design, for exponential family models, are given by the solutions to maximising a function $\phi(\bd)$ over $\mathcal{D}$. That is, a FIG-optimal design is given by $\tilde{\bD} = \left(\tilde{\bd}_1, \dots, \tilde{\bd}_n \right)^\T$ where $\tilde{\bd}_i$ maximises
$\phi(\bd)$ over $\mathcal{D}$, for $i=1,\dots,n$. Ostensibly, this is appealing since an $nk$-dimensional maximisation problem has been converted into a $k$-dimensional maximisation problem. However, the number of global maxima of $\phi(\bd)$ over $\mathcal{D}$ determines the maximum number of support points for a FIG-optimal design. If this number is less than $p$ (the number of parameters), then the design is under-supported, and we discuss how this is an undesirable feature for a default utility function. Note that for the remainder of this paper, the term maxima is used to mean global maxima.

\subsection{Form of $\phi(\bd)$}

\subsubsection{Known scale parameter}

We begin by considering the case where the scale parameter $\gamma$ is known. The following lemma provides the form for $\phi(\bd)$.

\begin{lemma} \label{theorem1}
For the exponential family models described in Section~\ref{sec:MODELS} with scale parameter $\gamma$ known, then
\begin{equation}
\phi(\bd) = \sum_{j=1}^p \mathrm{E}_{\btheta} \left[ \frac{1}{\mathrm{var}_{y \vert \btheta, \bd} (y)} \left( \frac{\partial \mu( \btheta, \bd)}{\partial \theta_j} \right)^2 \right].
\label{eqn:phid}
\end{equation}
The proof is given in Section~\ref{sec:suppproof1} of the Supplementary Material.
\end{lemma}

\subsubsection{Unknown scale parameter}

When the scale parameter $\gamma$ is unknown and of direct interest, we consider an extended $(p+1) \times 1$ parameter vector $\left(\btheta^{\T}, \gamma \right)^{\T}$. The following results gives the form of $\phi(\bd)$ in this case.

\begin{lemma} \label{theorem1b}
For the exponential family models described in Section~\ref{sec:MODELS} with scale parameter $\gamma$ unknown and of direct interest, then $\phi(\bd)$ is given by
\begin{equation}
\phi(\bd) = \sum_{j=1}^p \mathrm{E}_{\btheta, \gamma} \left[ \frac{1}{\mathrm{var}_{y \vert \btheta, \gamma, \bd} (y)} \left( \frac{\partial \mu( \btheta, \bd)}{\partial \theta_j} \right)^2 \right].
\label{eqn:phid1b}
\end{equation}
The proof is given in Section~\ref{sec:suppproof1b} of the Supplementary Material.
\end{lemma}

%When the scale parameter $\gamma$ is unknown and of direct interest, the element on the diagonal of the expanded $(p+1) \times (p+1)$ Fisher information matrix does not depend on the design \citep[see, for example,][]{kosmidis_etal_2020}, therefore
%$$\phi(\bd) = \sum_{j=1}^p \mathrm{E}_{\btheta, \gamma} \left[ \frac{1}{\mathrm{var}_{y \vert \btheta, \bd} (y)} \left( \frac{\partial \mu( \btheta, \bd)}{\partial \theta_j} \right)^2 \right].$$

When the scale parameter $\gamma$ is an unknown nuisance parameter, the situation is not as straightforward since the form of the marginal likelihood (\ref{eqn:marglik}) will depend on the form of $a(\cdot)$, $b(\cdot)$, $c_1(\cdot)$ and $c_2(\cdot)$. We focus on the case when the member of the exponential family is the normal distribution, with the variance of $y_i$ given by $\gamma$, for $i=1,\dots,n$. An inverse-gamma (IG) prior distribution for $\gamma$ is assumed giving a closed form expression for the marginal likelihood. If $\gamma \sim \mathrm{IG}(s_1/2, s_2/2)$ then the pdf is 
$$\pi(\gamma) = \frac{s_2^{\frac{s_1}{2}}}{2^{\frac{s_1}{2}}\Gamma\left(\frac{s_1}{2}\right)} \gamma^{-\frac{s_1}{2} - 1} \exp \left( - \frac{s_2}{2\gamma} \right)$$
for known constants $s_1$ (shape) and $s_2$ (scale). For such a model, the following lemma provides the form for $\phi(\bd)$.

\begin{lemma} \label{theorem2}
Under a normal distribution for the response with unknown scale parameter assumed to have prior distribution $\gamma \sim \mathrm{IG}(s_1/2,s_2/2)$, then  
\begin{equation}
\phi(\bd) = \sum_{j=1}^p \mathrm{E}_{\btheta} \left[ \left( \frac{\partial \mu(\btheta, \bd)}{\partial \theta_j} \right)^2 \right].
\label{eqn:phid2}
\end{equation}
The proof is given in Section~\ref{sec:suppproof2} of the Supplementary Material.
\end{lemma}

\subsection{Parameter redundancy in FIG-optimal designs}

A consequence of Lemmas~\ref{theorem1} to \ref{theorem2} is that to find the design points of any FIG-optimal design, we merely need to find vectors of controllable variables that maximise $\phi(\bd)$ over $\mathcal{D}$. Indeed, a FIG-optimal design exists with just one support point, i.e. $\tilde{\bD} = \left(\tilde{\bd}, \dots, \tilde{\bd}\right)^\T$ with $\tilde{\bd}$ any vector such that $\tilde{\bd} = \arg \max_{\bd \in \mathcal{D}} \phi(\bd)$. The number of maxima of $\phi(\bd)$ over $\mathcal{D}$ is the maximum number of support points $q$ of any FIG-optimal design. If $q < p$, then all FIG-optimal designs are under-supported, i.e. they have less support points than the number of unknown parameters. For such under-supported designs, the number of linearly independent rows of $M(\btheta, \bD)$ is at most $q<p$, meaning that $M(\btheta, \bD)$ is rank deficient, the model is parameter redundant and is therefore non-identifiable. 
%
%For such designs, the $n \times p$ matrix $M$, given by (\ref{eqn:Mmat}), will be rank-deficient, i.e. $\mathrm{rank}(M)<p$. \citet{catchpolemorgan1997} show, for exponential family models, if $M$ is rank-deficient, then the model is parameter redundant: unique classical estimates of $\btheta$ do not exist and the Fisher information matrix is singular for all $\btheta$. In the case of Bayesian inference, the issue is more nuanced. Under a proper prior distribution, the posterior distribution for $\btheta$ will exist for an under-supported design. However, parameter redundancy leads to Bayesian non-identifiability \citep[see, for example, ][]{gelfandsahu_1999}. This is where $\btheta$ can be decomposed as $\btheta = \left(\btheta_1^\T, \btheta_2^\T \right)^\T$ and the posterior distribution of $\btheta_1$, conditional on $\btheta_2$, is equal to the prior distribution of $\btheta_1$, conditional on $\btheta_2$. This does not mean that no learning has taken place about $\btheta_1$, rather that learning comes indirectly through $\btheta_2$.
%%\citep[see, for example, ][page 70]{OF2004}. 

In general, for an arbitrary statistical model (i.e. specification of response distribution and mean function $\mu(\btheta, \bd)$) and prior distribution for $\btheta$, analytically determining the number of maxima of $\phi(\bd)$, and therefore, whether or not all FIG-optimal designs are under-supported, will not be possible. In these cases, a numerical search of the maxima of $\phi(\bd)$ will be necessary. 

Note that even if a FIG-optimal design exists that is not under-supported, it can still result in parameter redundancy. This is because even though $\tilde{\bD}$ has at least $p$ support points, there still may be less than $p$ linearly independent rows of $M(\btheta, \bD)$. This issue is demonstrated in the context of (generalised) linear models in Sections~\ref{sec:normallinmod} and \ref{sec:Poisson}, where sufficient conditions are provided for parameter redundancy.

\section{Examples} \label{sec:EXAMPLES}

In this section we consider special cases of exponential family models such that Lemmas~\ref{theorem1} and~\ref{theorem2} can be applied, namely normal linear, Poisson, and logistic regression models and a non-linear compartmental model. These models have been frequently used in the literature to demonstrate developments of design of experiments methodology. In some cases, it is possible to determine conditions that result in parameter redundancy, and in others it is not (necessitating a numerical search).

\subsection{Normal linear regression model} \label{sec:normallinmod}

Design of experiments for the normal linear regression model is ubiquitous in the literature due to the applicability of such models for a wide range of applications and the availability of closed form results \citep[see, for example,][]{morris2011}. For such models, the exponential family distribution is normal and $\mu(\btheta, \bd) = f(\bd)^\T \btheta$ for a regression function $f : \mathcal{D} \to \mathbb{R}^p$. In this section, the elements of the regression function $f(\bd) = \left(f_1(\bd), \dots, f_p(\bd)\right)^\T$ are characterised as a monomial of the elements  of $\bd$, i.e.
$$f_j(\bd) = \prod_{r=1}^k d_{r}^{u_{rj}},$$
for $j=1,\dots,p$, where $u_{rj}$ is a non-negative integer. This general formulation allows incorporation of an intercept, main effects, interactions and powers of the controllable variables. The columns of the resulting matrix $[u]_{rj}$ should be unique. It is assumed that the row sums of $[u]_{rj}$ are positive which is equivalent to at least one element of $f(\bd)$ being dependent on each element of $\bd$. This is a weak assumption since it would be unwise to be able to control a variable but automatically exclude it from any terms in the model.

Suppose that the design space is constrained so that $d_{r} \in [a_r,b_r]$, with $a_r < b_r$, for $r=1,\dots,k$, and $\mathcal{D} = [a_1,b_1] \otimes [a_2,b_2] \otimes \dots \otimes [a_k, b_k]$ where $\otimes$ denotes the Cartesian product. 

As an example, consider the second order response surface model in $k=2$ controllable variables in \citet[pages 170-171]{ADT2007}. Here $f(\bd) = \left(1, d_{1},d_{2},d_{1}^2, d_{2}^2, d_{1}d_{2}\right)^\T$, corresponding to a model with an intercept, main effects and squared terms for both controllable variables and an interaction, with $p=6$. In this case
$$\left[u\right]_{rj} = \left( \begin{array}{cccccc}
0 & 1 & 0 & 2 & 0 & 1\\
0 & 0 & 1 & 0 & 2 & 1 \end{array}\right).$$
The design space is constrained such that $a_1=a_2 = -1$ and $b_1=b_2 =1$, with $\mathcal{D} = [-1,1]^2$.

Under the above setup, FIG-optimal designs can result in parameter redundancy and, therefore, non-identifiability. To see this, first consider the following result.

\begin{proposition} \label{lemma3}
Under the normal linear model described above with an assumed inverse-gamma prior distribution for $\gamma$, any FIG-optimal design is $\tilde{\bD} = \left(\tilde{\bd}_1,\dots,\tilde{\bd}_n\right)^\T$ where $\tilde{\bd}_i = \left(\tilde{d}_{i1},\dots, \tilde{d}_{ik}\right)^\T$ with
\begin{equation}
\tilde{d}_{ir} = \left\{
\begin{array}{ll}
a_r & \mbox{if $|a_r| > |b_r|$;}\\
b_r & \mbox{if $|a_r| < |b_r|$;}\\
\pm |b_r| & \mbox{if $|a_r| = |b_r|$;}
\end{array} \right.
\label{eqn:lemma3}
\end{equation}
for $i=1,\dots,n$ and $r=1,\dots,k$.
The proof is given in Section~\ref{sec:suppproof3} of the Supplementary Material.
\end{proposition}

A consequence of Proposition~\ref{lemma3} is that the maximum number of support points for any FIG-optimal design is $q = 2^C$ where $C = \sum_{r=1}^k I(\vert a_r \vert = \vert b_r \vert)$ with $I(A)$ denoting the indicator function of event $A$. Therefore, all FIG-optimal designs are under-supported if $p > 2^C$. Commonly, $a_r = -b_r$, for $r=1,\dots,k$ and all FIG-optimal designs are under-supported if $p > 2^k$. 

Consider the second-order response model example above with $k=2$ controllable variables. Since $6 = p > 2^k = q= 4$, all FIG-optimal designs for this model are under-supported resulting in parameter redundancy. Specifically, the design points are $\tilde{\bd}_i = \left( \pm 1, \pm 1\right)^\T$.

Even if $q \ge p$, and therefore a FIG-optimal design exists which is not under-supported, all FIG-optimal designs can still result in parameter redundancy. The following result provides a sufficient condition on the form of the model which results in this outcome.

\begin{proposition} \label{lemma4}
Under the normal linear model described above with an assumed inverse-gamma prior distribution for $\gamma$, suppose there exists $j_1,j_2 \in \left\{ 1, \dots, k \right\}$, with $j_1 \ne j_2$, such that for all $r=1,\dots,k$, $u_{rj_1} = 0$, and there exists an $\bar{r} \in \left\{1, \dots,k\right\}$ such that $u_{\bar{r}j_2} =2$ and $u_{rj_2}=0$ for $r \in \left\{ 1, \dots, k \vert r \ne \bar{r} \right\}$. Then the model will be parameter redundant. The proof is given in Section~\ref{sec:suppproof4} of the Supplementary Material.
\end{proposition}

Proposition~\ref{lemma4} means that if a normal linear regression model includes an intercept and a term given by the square of any of the controllable variables, then any FIG-optimal design results in a parameter redundant model, irrespective of whether or not it is under-supported. An intuitive explanation for this result is as follows. From Proposition~\ref{lemma3}, it can be seen that the design points of a FIG optimal design are at the extremes of the design space. When a controllable variable is constrained to be in an interval symmetric about zero, i.e. $|a_r| = |b_r|$, to be able to uniquely estimate the parameters of a model which includes a squared term for this controllable variable (and an unknown intercept), design points in the interior of the interval are required to be able to estimate curvature.

As a particular example, suppose there is $k=1$ controllable variable $d \in [-1,1]$ and the model has
$$f(d) = \left(1, d \right)^{\T} \qquad \mbox{} \qquad \left[u\right]_{rj} = \left(\begin{array}{ll} 0 & 2 \end{array} \right),$$
with $p=2$. All FIG-optimal designs have the form $\tilde{\bD} = \left\{ \tilde{d}_1, \dots, \tilde{d}_n \right\}$ with $\tilde{d}_i = \pm 1$. Now suppose that at least one design point is $-1$ and at least one design point is $+1$. Then any design of this form is FIG-optimal but not under-supported (since it has two support points). However, by Proposition~\ref{lemma4}, the model is parameter redundant. To see this, note that the two columns of the $n \times 2$ matrix $M(\btheta, \bD)$ are both $\bone_n$ (an $n \times 1$ vector of ones), so the matrix is rank deficient. 

%More intuitively, the classical maximum likelihood estimates of $\btheta$ are given by $\left(M(\btheta, \bD)^{\T}M(\btheta, \bD)\right)^{-1}M(\btheta, \bD)^{\T}\by$. However $M(\btheta, \bD)^{\T}M(\btheta, \bD)$ is singular due to the rank deficiency of $M(\btheta, \bD)$.

%Furthermore, in general, all FIG-optimal designs will lead to parameter redundancy if there exist $j_1,j_2 = 1, \dots p$ with $j_1 \ne j_2$ where (i) $u_{rj_1} = 0$ and $u_{rj_2} >1$ is even and/or (ii) $u_{rj_1} = 1$ and $u_{rj_2} >1$ is odd. In other words, for models with an intercept and a main effect (first-order term) for each controllable variable, inclusion of higher-order powers of the controllable variables will lead to parameter redundancy.
%
%Consider a second-order response surface model in $k$ controllable variables where $a_r = -1$ and $b_r=1$ for $r=1,\dots,k$. For a model with an intercept, main effects and squared terms for all controllable variables, and all two-way interactions, there are $p = 1 + 2k + {k \choose 2}$ parameters. For $k \le 4$, $q = 2^k >p$ and a FIG-optimal design will exist that is not under-supported. However, due to the inclusion of an intercept and squared terms, the model is parameter redundant.

\subsection{Poisson regression model} \label{sec:Poisson}

Under the Poisson generalised linear model, the exponential family distribution is Poisson and the mean function is $\mu(\btheta, \bd) = \exp \left( f(\bd)^\T \btheta \right)$, where $f(\bd)$ is a regression function as defined in Section~\ref{sec:normallinmod}. The scale parameter is $\gamma = 1$ (i.e. known) and $\mathrm{var}(y \vert \btheta, \bd) = \mu(\btheta, \bd)$. From Lemma~\ref{theorem1}, 
$$\phi(\bd) = \mathrm{E}_{\btheta} \left[ \exp \left(f(\bd)^\T \btheta \right) \right] \sum_{j=1}^p f_j(\bd)^2.$$
For an arbitrary prior distribution it will not be possible to determine the number of maxima of $\phi(\bd)$. 

Instead, we assume that the elements of $\btheta$ are, a-priori, independent with $\theta_j \sim \mathrm{N}\left(0, \sigma^2_j\right)$. In this case,
$$\phi(\bd) = \exp \left( \frac{1}{2} \sum_{j=1}^p \sigma_j^2 f_j(\bd)^2\right)  \sum_{j=1}^p f_j(\bd)^2.$$
If the design space is constrained so that $d_{ir} \in [a_r,b_r]$, for $i=1,\dots,n$, then using a similar argument to the proof of Proposition~\ref{lemma3}, a FIG-optimal design is $\tilde{\bD} = \left(\tilde{\bd}_1,\dots,\tilde{\bd}_n\right)^\T$ where $\tilde{\bd}_i = \left(\tilde{d}_{i1},\dots, \tilde{d}_{ik}\right)^\T$ and $\tilde{d}_{ir}$ is given by (\ref{eqn:lemma3}), for $i=1,\dots,n$ and $r=1,\dots,k$. The consequences of using such a design for a normal linear regression model apply here, i.e. parameter redundancy if i) the design is under-supported when $2^C < p$ (see Proposition~\ref{lemma3}), and/or ii) inclusion of an intercept and squared term (see Proposition~\ref{lemma4}).

\subsection{Logistic regression model} \label{sec:logreg}
Consider a logistic regression model where $y_i \sim \mathrm{Bernoulli}\left(\mu(\btheta, \bd_i)\right)$, for $i=1,\dots,n$ where $\mu(\btheta, \bd) = 1/\left( 1 + \exp \left( - f(\bd)^\T \btheta \right) \right)$. The scale parameter is $\gamma = 1$ and $\mathrm{var}(y \vert \btheta, \bd) = \mu(\btheta, \bd) \left(1 - \mu(\btheta, \bd) \right)$. From Lemma~\ref{theorem1},
\begin{equation}
\phi(\bd) = \mathrm{E}_{\btheta} \left[ \mu(\btheta, \bd) \left(1 - \mu(\btheta, \bd) \right) \right] \sum_{j=1}^p f_j(\bd)^2.
\label{eqn:logregphi}
\end{equation}

Similar to the Poisson regression model in Section~\ref{sec:Poisson}, under an arbitrary prior distribution for $\btheta$, it will not be possible to evaluate $\phi(\bd)$ in closed form and thus determine the number of maxima. 

Instead we consider a particular example from \citet{W2006} which is frequently used to demonstrate design of experiments methodology \citep[see, for example,][]{GJS2009}. In this example, there are $k=4$ controllable variables $\bd = \left(d_1,\dots,d_4\right)^\T$. The regression function is $f(\bd) = \left(1, d_1,\dots, d_4\right)^\T$ leading to an intercept, a main effect for each controllable variable and $p=5$ parameters. Independent uniform prior distributions are assumed for the elements of $\btheta$ as follows
$$
\begin{array}{cccccccc}
\theta_1 \sim \mathrm{U}[-3,3]\,, & &\theta_2\sim\mathrm{U}[4,10]\,, & &\theta_3\sim\mathrm{U}[5,11]\,,\\ 
\theta_4\sim \mathrm{U}[-6,0]\,, & & \theta_5 \sim \mathrm{U}[-2.5,3.5]\,.\\ 
\end{array}
$$
The elements of $\bd$ are constrained such that $a_r = -1$ and $b_r = 1$, for $r=1,\dots,4$. 

The prior expectation involved in evaluating the function $\phi(\bd)$ in (\ref{eqn:logregphi}) is approximated using Gauss-Legendre quadrature \citep[see, for example,][Chapter 3]{DR1984}, i.e. 
$$\hat{\phi}(\bd) = \left[\sum_{c=1}^C \omega_c \mu(\bt_c, \bd) \left(1 - \mu(\bt_c, \bd) \right)\right] \sum_{j=1}^p f_j(\bd)^2,$$
where $\left\{ \bt_1,\dots,\bt_C \right\}$ are quadrature points with weights $\left\{\omega_1, \dots, \omega_C \right\}$. The number of quadrature points is chosen to be $C=10,000$.

To determine the number of maxima of $\phi(\bd)$, the function $\hat{\phi}(\bd)$ is numerically maximised using a quasi-Newton scheme, independently from $B=1000$ different starting points, generated uniformly over $\mathcal{D}$. From this procedure, it is found that $\hat{\phi}(\bd)$ has only two maxima, at $\tilde{d} = \pm \left( -0.8701,1,1,1\right)^{\T}$ (to 4 decimal places). That means that any FIG-optimal design for the above logistic regression model has a maximum of $q=2$ support points, thus resulting in parameter redundancy.

\subsection{Compartmental non-linear model} \label{sec:como}

Compartmental models are frequently used in pharmacokinetics to model the uptake and excretion of a substance from an organism. Here the first-order compartmental model from \citet{GJS2009} is considered. There is $k=1$ controllable variable denoted $d \in \mathcal{D}=[0,24]$, representing the time (in hours) from a substance first entering the organism until the response ($y$; the concentration of substance in the organism) is measured. The response is assumed to be normally distributed with mean function $\mu(\btheta, d) = \theta_3 \left( \exp (- \theta_1 d) - \exp (- \theta_2 d) \right)$, where $\theta_2 >\theta_1$. Following \citet{GJS2009}, a-priori the elements of $\btheta$ are assumed independent with $\theta_1 \sim U(0.01884, 0.09884)$, $\theta_2 \sim U(0.298, 8.298)$, and $\theta_3$ having a point mass at 21.8. Finally, the scale parameter $\gamma$ (the common variance of the response) is assumed to have an inverse gamma prior distribution.

%\begin{center}
%\begin{figure}[t!]
%\centerline{\includegraphics[scale=0.6]{compplot}}
%\caption{Plot of $\phi(d)$ against $d$ for the compartmental model.} \label{fig:compplot}
%\end{figure}
%\end{center}

From Lemma~\ref{theorem2},
$$
\phi(d) = \mathrm{E}_{\btheta} \left[ \sum_{j=1}^p \left(\frac{\partial \mu(\btheta, d)}{\partial \theta_j}\right)^2 \right].
$$
Similar to Section~\ref{sec:logreg}, $\phi(d)$ is approximated using Gauss-Legendre quadrature. Figure~\ref{fig:compplot} in the Supplementary Material shows the resulting approximation $\hat{\phi}(d)$ plotted against $d \in [0,24]$. From inspection, $\hat{\phi}(d)$ has one maximum at $\tilde{d} = 24$. Therefore for a design of $n$ sampling times, the unique FIG-optimal design is $\tilde{D} = \tilde{d} \bone_n$. This design is under-supported leading to parameter redundancy.

\section{Concluding remarks} \label{sec:DISC}

In this paper it is shown that, for exponential family models the design points of FIG-optimal designs can be found by maximising a function $\phi(\bd)$ over $\mathcal{D}$. The advantage is that this significantly simplifies the identification of a FIG-optimal design. However, the disadvantage is that if the number of maxima of $\phi(\bd)$ is less than the number of parameters then all FIG-optimal designs will be under-supported resulting in parameter redundancy and, therefore, a non-identifiable model. Furthermore, even when a FIG-optimal design exists that is not under-supported, it can still result in parameter redundancy for (generalised) linear models on the inclusion of higher-order terms.

The conclusions from the paper are as follows. 
\begin{itemize}
\item
If the experimental aim is accurately represented by Fisher information gain, and a Bayesian analysis pursued post-experiment, then the findings are inconsequential. Under a proper prior distribution, the posterior distribution for $\btheta$ will exist for a parameter-redundant model and a FIG-optimal design maximises expected gain in information as measured by the FIG utility. 
\item
However, if a default utility is to be considered, with the possibility of a non-Bayesian analysis post-experiment, then Fisher information gain should be used with caution.
\end{itemize}

In the latter case, it should be verified that the FIG-optimal design does not result in parameter redundancy. Alternatively, a different approach to design of experiments whilst still incorporating Fisher information gain could be employed. \cite{prangle_etal_2020} empirically observed some of the disadvantages of FIG-optimal designs explored in this paper. To counter this, \cite{prangle_etal_2020} developed a game theoretic approach to design of experiments which is still based on Fisher information gain, but circumvents the disadvantages.

\setcounter{section}{0}
\setcounter{equation}{0}
\def\theequation{SM\arabic{section}.\arabic{equation}}
\def\thesection{SM\arabic{section}}
\def\thefigure{SM\arabic{figure}}

\section*{Supplementary Material}

\section{Proof of Lemma~\ref{theorem1}} \label{sec:suppproof1}
%The expected FIG utility is
%\begin{eqnarray*}
%U_{FIG}(\bD) & = & \mathrm{E}_{\btheta} \left[ \mathrm{tr} \left\{ \mathcal{I}(\btheta , \bD) \right\} \right]\\
%& = & \mathrm{E}_{\btheta} \left[ \mathrm{tr} \left\{ \mathrm{var}_{\by \vert \btheta, \bD} \left( \frac{\partial \log \pi(y_i \vert \btheta, \bd_i)}{\partial \btheta} \right) \right\} \right]\\
%& = & \sum_{i=1}^n \sum_{j=1}^p \mathrm{E}_{\btheta} \left[ \mathrm{var}_{y_i \vert \btheta, \bd_i} \left( \frac{\partial \log \pi(y_i \vert \btheta, \bd_i)}{\partial \theta_j} \right) \right].
%\end{eqnarray*}
%The last line follows from the assumption of the elements of $\by$ under exponential family models. Due to the fungible nature (see Section~\ref{sec:MODELS}) of the elements of $\by$,
%$$U_{FIG}(\bD) =  \sum_{i=1}^n \phi(\bd_i),$$
%where
%\begin{eqnarray*}
%\phi(\bd) & = & \sum_{j=1}^p \mathrm{E}_{\btheta} \left[ \mathrm{var}_{y \vert \btheta, \bd} \left( \frac{\partial \log \pi(y \vert \btheta, \bd_i)}{\partial \theta_j} \right) \right]\\
%& = & \sum_{j=1}^p \mathrm{E}_{\btheta} \left[ \frac{1}{\mathrm{var}_{y \vert \btheta, \bd} (y)} \left( \frac{\partial \mu(\btheta, \bd)}{\partial \theta_j} \right)^2 \right].
%\end{eqnarray*}

The expected FIG utility is 
\begin{eqnarray*}
U_{FIG}(\bD) & = & \mathrm{E}_{\btheta} \left[ \mathrm{tr} \left\{ \mathcal{I}(\btheta , \bD) \right\} \right]\\
& = & \mathrm{E}_{\btheta} \left[ \sum_{j=1}^p \sum_{i=1}^n \frac{1}{\mathrm{var}_{y_i \vert \btheta, \gamma, \bd_i} \left(y_i\right)} \left( \frac{\partial \mu(\btheta, \bd_i)}{\partial \theta_j}\right)^2 \right]
\end{eqnarray*}
where the second line follows from the assumption that the elements of $\by$ follow an exponential family distribution. Interchanging the order of the summations and the expectation, 
\begin{eqnarray*}
U_{FIG}(\bD) & = & \sum_{i=1}^n \sum_{j=1}^p\mathrm{E}_{\btheta} \left[\frac{1}{\mathrm{var}_{y_i \vert \btheta, \gamma, \bd_i} \left(y_i\right)} \left( \frac{\partial \mu(\btheta, \bd_i)}{\partial \theta_j}\right)^2 \right]\\
& = & \sum_{i=1}^n \phi(\bd_i),
\end{eqnarray*}
with $\phi(\bd)$ given by (\ref{eqn:phid}). 

%Due to the fungible nature of the elements of $\by$, to maximise $U_{FIG}(\bD)$, it is equivalent to maximise each $\phi(\bd_i)$.

\section{Proof of Lemma~\ref{theorem1b}} \label{sec:suppproof1b}

The Fisher information matrix for the extended parameter vector $\left(\btheta^{\T}, \gamma \right)^{\T}$ is a $(p+1) \times (p+1)$ matrix. The $p \times p$ sub-matrix corresponding to $\btheta$ is identical to the case when $\gamma$ is known. The element on the diagonal of the Fisher information matrix corresponding to $\gamma$ is $na''(-1/\gamma)/2\gamma^4$ \citep[see, for example,][]{kosmidis_etal_2020}. Therefore,
$$U_{FIG}(\bD) = \mathrm{E}_{\btheta, \gamma} \left[ \sum_{j=1}^p \sum_{i=1}^n \frac{1}{\mathrm{var}_{y_i \vert \btheta, \gamma, \bd_i} \left(y_i\right)} \left( \frac{\partial \mu(\btheta, \bd_i)}{\partial \theta_j}\right)^2 \right] + \mathrm{E}_{\btheta, \gamma} \left[ n\frac{a''\left(-\frac{1}{\gamma}\right)}{2\gamma^4}\right].$$
In the first term, the order of the summations and the expectation are interchanged. The second term does not depend on the design $\bD$. Therefore, up to an additive constant not depending on the design,
\begin{eqnarray*}
U_{FIG}(\bD) & = & \sum_{i=1}^n \sum_{j=1}^p\mathrm{E}_{\btheta, \gamma} \left[\frac{1}{\mathrm{var}_{y_i \vert \btheta, \gamma, \bd_i} \left(y_i\right)} \left( \frac{\partial \mu(\btheta, \bd_i)}{\partial \theta_j}\right)^2 \right]\\
& = & \sum_{i=1}^n \phi(\bd_i),
\end{eqnarray*}
with $\phi(\bd)$ given by (\ref{eqn:phid1b}). 

%Due to the fungible nature of the elements of $\by$, to maximise $U_{FIG}(\bD)$, it is equivalent to maximise each $\phi(\bd_i)$.

\section{Proof of Lemma~\ref{theorem2}} \label{sec:suppproof2}

The marginal likelihood (\ref{eqn:marglik}) is 
\begin{eqnarray}
\pi(\by \vert \btheta, \bD) & = & \int \pi(\by \vert \btheta, \gamma, \bD) \pi(\gamma) \mathrm{d}\gamma \nonumber \\
& \propto & \int \gamma^{-n/2} \exp \left[ - \frac{\sum_{i=1}^n \left(y_i - \mu(\btheta, \bd_i)\right)^2}{2\gamma} \right] \gamma^{-s_1/2-1} \exp \left( - \frac{s_2}{2\gamma} \right) \mathrm{d}\gamma \nonumber \\
& & \label{eqn:bi1}\\
& = & \frac{2^{\frac{s_1+n}{2}} \Gamma \left(\frac{s_1+n}{2}\right)}{\left[ s_2 + \sum_{i=1}^n (y_i - \mu(\btheta, \bd_i))^2 \right]^\frac{s_1+n}{2}} \label{eqn:bi2} \\
& \propto & \left[ 1 + \frac{1}{s_2} \left(\by - \bmu(\btheta, \bD)\right)^{\T}\left(\by - \bmu(\btheta, \bD)\right) \right]^{-\frac{s_1+n}{2}}, \label{eqn:bi3}
\end{eqnarray}
where $\bmu(\btheta, \bD) = \left(\mu(\btheta, \bd_1), \dots, \mu(\btheta, \bd_n)\right)^{\T}$. The expression in (\ref{eqn:bi2}) follows since the integrand in (\ref{eqn:bi1}) is proportional to the pdf of an inverse-gamma distribution with shape parameter $s_1+n$ and scale parameter $s_2 +  \sum_{i=1}^n (y_i - \mu(\btheta, \bd_i))^2$. The final expression in (\ref{eqn:bi3}) is the pdf of the multivariate t-distribution \citep[see, for example, ][page 1]{KN2004}, with location vector $\bmu(\btheta, \bD)$, scale matrix $(s_2/s_1)I_n$ and $s_1$ degrees of freedom. From \cite{lange_etal_1989} the Fisher information under this marginal model is 
\begin{equation}
\mathcal{I}(\btheta, \bD) = \frac{s_1(s_1+n)}{s_2(s_1+n+2)} M(\btheta, \bD)^\T M(\btheta, \bD).
\label{eqn:tfim}
\end{equation}
Using the same argument as in the proof of Lemmas~\ref{theorem1} and~\ref{theorem1b}, it follows that the expected FIG utility is
$$U_{FIG}(\bD) = \frac{s_1(s_1+n)}{s_2(s_1+n+2)} \sum_{i=1}^n \phi(\bd_i),$$
with $\phi(\bd)$ given by (\ref{eqn:phid2}) in the main manuscript. 

%Due to the fungible nature of the elements of $\by$, to maximise $U_{FIG}(\bD)$, it is equivalent to maximise each $\phi(\bd_i)$.

 \begin{center}
\begin{figure}[t!]
\centerline{\includegraphics[scale=0.6]{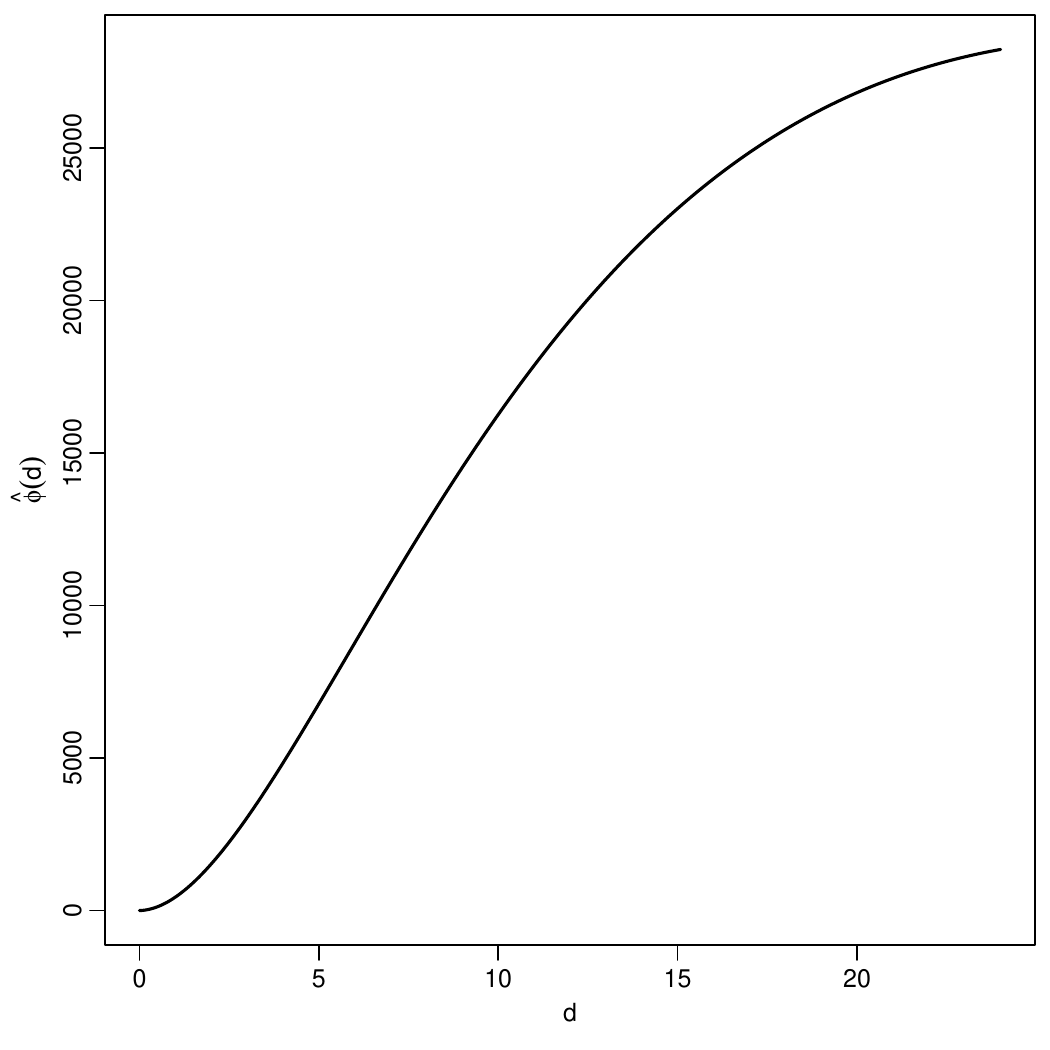}}
\caption{Plot of $\hat{\phi}(d)$ against $d$ for the compartmental non-linear model.} \label{fig:compplot}
\end{figure}
\end{center}

\section{Proof of Proposition~\ref{lemma3}} \label{sec:suppproof3}

Since $\partial \mu(\btheta, \bd)/\partial \theta_j = f_j(\bd)$, independent of $\btheta$, it follows from Lemma~\ref{theorem2} that
$$\phi(\bd) = \sum_{j=1}^p \prod_{r=1}^k \left(d_{r}^{u_{rj}}\right)^2.$$
The positive row sum assumption, along with with $u_{rj} \ge 0$, means that $\phi(\bd)$ is maximised by 
$$\tilde{d}_{r} = \left\{
\begin{array}{ll}
a_r & \mbox{if $|a_r| > |b_r|$;}\\
b_r & \mbox{if $|a_r| < |b_r|$;}\\
\pm |b_r| & \mbox{if $|a_r| = |b_r|$;}
\end{array} \right.$$
for $r=1,\dots,k$.

\section{Proof of Proposition~\ref{lemma4}} \label{sec:suppproof4}

Under a FIG optimal design, the $n \times p$ matrix $M(\btheta,\tilde{\bD})$ has $ij$th element given by $f_j(\tilde{\bd}_i)$. From Proposition~\ref{lemma3}, the $j_1$th column of $M(\btheta,\tilde{\bD})$ is $\bone_n$ (where $\bone_n$ denotes the $n \times 1$ vector of ones) and the $j_2$th column is 
$$\left(f_{j_2}(\tilde{\bd}_1), \dots, f_{j_2}(\tilde{\bd}_n)\right)^{\T} = \left\{ \begin{array}{ll}
a_{\bar{r}}^2 \bone_n & \mbox{if $|a_{\bar{r}}| > |b_{\bar{r}}|$;}\\
b_{\bar{r}}^2 \bone_n & \mbox{if $|a_{\bar{r}}| < |b_{\bar{r}}|$;}\\
b_{\bar{r}}^2 \bone_n & \mbox{if $|a_{\bar{r}}| = |b_{\bar{r}}|$.}
\end{array} \right.$$
In all cases, columns $j_1$ and $j_2$ are linearly dependent. Therefore $M(\btheta,\tilde{\bD})$ has less than $p$ linearly independent columns and is rank deficient.

\section{Plot of $\hat{\phi}(d)$ against $d$ for the compartmental non-linear model}

Figure~\ref{fig:compplot} shows the approximate objective function $\hat{\phi}(d)$, for the compartmental non-linear model, plotted against $d \in [0,24]$.

\bibliographystyle{rss}      % Chicago style, author-year citations
\bibliography{mybib}   % name your BibTeX data base

\end{document}